\documentclass{amsart}
\usepackage{amssymb}
\usepackage{amsmath}

\allowdisplaybreaks

\theoremstyle{plain}
\newtheorem{theorem}{Theorem}

\newtheorem{lemma}[theorem]{Lemma}

\theoremstyle{definition}

\theoremstyle{remark}

\numberwithin{equation}{section}
\numberwithin{theorem}{section}

\newcommand{\qbin}[3]{\genfrac{[}{]}{0pt}{}{#1}{#2}_{#3}}

\setbox0=\hbox{$+$}

\newdimen\plusheight
\plusheight=\ht0
\def\+{\;\lower\plusheight\hbox{$+$}\;}

\setbox0=\hbox{$-$}
\newdimen\minusheight
\minusheight=\ht0
\def\-{\;\lower\minusheight\hbox{$-$}\;}

\setbox0=\hbox{$\cdots$}
\newdimen\cdotsheight
\cdotsheight=\plusheight
\def\cds{\lower\cdotsheight\hbox{$\cdots$}}
\begin{document}
\title[ Some more identities of the Rogers-Ramanujan type ]
       { Some more identities of the Rogers-Ramanujan type}
\author{Douglas Bowman}
\address{ Northern Illinois University\\
   Mathematical Sciences\\
   DeKalb, IL 60115-2888
} \email{bowman@math.niu.edu}

\author{James Mc Laughlin}
\address{Mathematics Department\\
 Anderson Hall\\
West Chester University, West Chester, PA 19383}
\email{jmclaughl@wcupa.edu}

\author{Andrew V. Sills}
\address{Department of Mathematical Sciences\\
Georgia Southern University\\
Statesboro, GA 30460-8093}
\email{asills@georgiasouthern.edu}

\thanks{The research of the first author was partially supported by National Science Foundation grant
DMS-0300126.}

\date{\today}

\begin{abstract}
In this we paper we prove several new identities of the
Rogers-Ramanujan-Slater type. These identities were found as the
result of  computer searches. The proofs involve a variety of
techniques, including series-series identities, Bailey pairs, a
theorem of Watson on basic hypergeometric series, generating
functions and miscellaneous methods.
\end{abstract}

\maketitle


\section{introduction}
The most famous of the ``$q$-series=product" identities are the
Rogers-Ramanujan identities:
\begin{align}\label{rr1}
\sum_{n=0}^{\infty}\frac{q^{n^{2}}}{(q;q)_{n}}&=
\prod_{j=0}^{\infty}\frac{1}{(1-q^{5j+1})(1-q^{5j+4})},\\
\sum_{n=0}^{\infty}\frac{q^{n^{2}+n}}{(q;q)_{n}}&=
\prod_{j=0}^{\infty}\frac{1}{(1-q^{5j+2})(1-q^{5j+3})}. \label{rr2}
\end{align}
These identities have a curious history (\cite{H37}, p. 28). They
were first proved by L.J. Rogers in 1894 (\cite{R94}) in a paper
that was completely ignored. They were rediscovered (without proof)
by Ramanujan sometime before 1913. In 1917, Ramanujan rediscovered
Rogers's paper. Also in 1917, these identities were rediscovered and
proved independently by Issai Schur (\cite{S17}). They were also
discovered independently by R. Baxter  (see \cite{AB89} for
details). An account of the many proofs of the Rogers-Ramanujan
identities can be found in \cite{A89}.

There are numerous identities that are similar to the
Rogers-Ramanujan identities. These include identities by Jackson
(\cite{J28}), Rogers~(\cite{R94} and \cite{R17}) and
Bailey~(\cite{B47} and \cite{B49}). Of special note is Slater's 1952
paper \cite{S52}, which contains a list of 130 such identities, many
of them new (see the paper by the third author~\cite{S03}, for an
annotated version of Slater's list). There are also other identities
of Rogers-Ramanujan type in the literature.

In the present paper we describe the results of some numerical
investigations, which were undertaken with the aim of finding new
Rogers-Ramanujan type identities. One reason for searching for such
identities is the possibility that new identities might be found
experimentally which could \emph{not} be proven using presently
existing methods, thus necessitating developments in the
mathematical theory.

These investigations did uncover several new identities, although
all were provable using known methods. The proof of these identities
lets us present examples of the various methods used to prove
identities of the Rogers-Ramanujan type. These methods include using
series-series identities,  Bailey pairs, generating functions and
some miscellaneous methods.

One search involved computing to high precision series of the form
\begin{equation}\label{Seq}
S:=\sum_{n=0}^{\infty}\frac{q^{(a n^2+b
n)/2}(-1)^{cn}(d,e;q)_n}{(f,g,q;q)_n},
\end{equation}
for a fixed numerical value of $q$, for
\[
d,e,f,g \in \{0,-1,q,-q,-q^2,q^2\}, \hspace{20pt} c \in \{0,1\},
\]
and for integers $a$ and $b$. This choice for the form of $S$ was
motivated by the fact that many series on the Slater list have this
form.

For each particular choice of the parameters $a$, $b$, $c$, $d$,
$e$, $f$ and $g$, a numerical comparison was performed to see if
\begin{equation}\label{SPeq}
S-\prod_{j=1}^{L}(q^j;q^{L})_{\infty}^{s_j}=0,
\end{equation}
for integers $s_j$ and $L \in \{20,24,28,32,36\}$. With sufficient
precision, a small numerical value for the left side of \eqref{SPeq}
indicated either  a known identity or a potential new identity,
which then needed to be proved.

 We also tried searches where the series had the
forms
\begin{align*}
S'&:=\sum_{n=0}^{\infty}\frac{q^{(a n^2+b
n)/2}(-1)^{cn}(d,e;q^2)_n}{(f,g,q^2;q^2)_n},\\
S''&:=\sum_{n=0}^{\infty}\frac{q^{(a n^2+b
n)/2}(-1)^{cn}(d;q)_n}{(e;q^2)_{n+1}(q;q)_{n+1}},
\end{align*}
the latter form being motivated by identity (56) on Slater's list.

It is possible that other choices for the form of the series, or
extending the choices for the various parameters, may turn up new
identities. The computations associated with the searches described
above were performed using \emph{PARI/GP}.

\section{Some Identities Discovered using PARI/GP}
In this section we list the new identities found during the
\emph{PARI/GP} searches. We organize them according to the methods
used in their proofs.

\subsection{Infinite Series Transformations}
Infinite series transformations can be \, used to derive new
identities from known identities, since if one of the series can be
expressed as an infinite product, for certain values of the
parameters, then the other series automatically also has an
expression as an infinite product. Before coming to the identities
in this subsection, we list some infinite series transformations
that are necessary for the proofs.

We first recall  Heine's $q$-Gauss sum.
\begin{equation*}
\sum_{n=0}^{\infty}\frac{(a,b;q)_n}{(c,q;q)_{n}}\left
(\frac{c}{ab}\right)^n
=\frac{(c/a,c/b;q)_{\infty}}{(c,c/ab;q)_{\infty}}
\end{equation*}
 Let $b \to \infty$ to get
\begin{equation}\label{eq2}
\sum_{n=0}^{\infty} \frac{(a;q)_{n}q^{n(n-1)/2}(-c/a)^{n}}
{(c;q)_{n}(q;q)_{n}} =\frac{(c/a;q)_{\infty}}{(c;q)_{\infty}}.
\end{equation}

The second identity we need is the following.
\begin{equation}\label{eq1}
\sum_{n=0}^{\infty} \frac{(a;q)_{n}q^{n(n-1)/2}\gamma^{n}}
{(b;q)_{n}(q;q)_{n}} = \frac{(-\gamma;
q)_{\infty}}{(b;q)_{\infty}}\sum_{n=0}^{\infty} \frac{(-a \gamma
/b;q)_{n}q^{n(n-1)/2}(-b)^{n}} {(-\gamma;q)_{n}(q;q)_{n}},
\end{equation}
This identity, in a more symmetric form,  is found in Ramanujan's
lost notebook \cite{S88} and a proof can be found in the recent book
by Andrews and Berndt \cite{AB05}. It also follows from the second
iteration of Heine's transformation for a $_2 \phi_1$ series
\cite[equations III.1 and III.2, page 359]{GR04}. The form at
\eqref{eq1} is better suited to our present requirements than
Ramanujan's more symmetric form.

We also recall the following transformation \cite[page 80]{GR04}:
\begin{multline*}
\sum_{n=0}^{\infty}\frac{(1-aq^{4n})(a,b,c,cq,d,dq;q^2)_n}
{(1-a)(aq^2/b,aq^2/c,aq/c,aq^2/d,aq/d,q^2;q^2)_n}\left(\frac{a^2q^2}{bc^2d^2}\right)^n\\
=\frac{(aq,aq/bc;q)_{\infty}(a^2q^2/c^2d^2,a^2q^2/d^2b;q^2)_{\infty}}{
(aq/b,aq/c;q)_{\infty}(a^2q^2/d^2,a^2q^2/c^2d^2b;q^2)_{\infty}}\\
\times
\sum_{n=0}^{\infty}\frac{(1+a/dq^{2n})(-a/d,c;q)_n(b,aq/d^2;q^2)_n}
{(1+a/d)(-aq/cd,q;q)n(a^2q^2/d^2b,aq;q^2)_n}\left(\frac{aq}{bc}\right)^n.
\end{multline*}
Upon replacing $c$ with $aq/c$, then letting $a\to 0$ and finally
$b\to \infty$, we get
\begin{equation}\label{eq16}
\sum_{n=0}^{\infty} \frac{(d;q)_{2n}q^{n^2-n}(-c^2/d^2)^{n}}
{(q^2;q^2)_{n}(c;q)_{2n}}
 =
\frac{(c^2/d^2;q^2)_{\infty}}{( c;q)_{\infty}} \sum_{n=0}^{\infty}
\frac{q^{n^2-n}(-c)^{n}} {(q;q)_{n}(-c/d;q)_{n}}.
\end{equation}

\begin{theorem}\label{c2}
\begin{equation}\label{sereq1}
\sum_{n=0}^{\infty} \frac{(-q^2;q^2)_{n}q^{n^2}}{(q;q)_{2n+1}}
=\frac{(-q;q^2)_{\infty}}{(q;q^2)_{\infty}}.
\end{equation}
\begin{equation}\label{sereq2}
\sum_{n=0}^{\infty}
\frac{(-1;q)_{2n}q^{n^2+n}}{(q^2;q^2)_{n}(q^2;q^4)_{n}}
=\frac{(-q^3;q^6)_{\infty}^2(q^6;q^6)_{\infty}(-q^2;q^2)_{\infty}}{(q^2;q^2)_{\infty}}.
\end{equation}
\begin{equation}\label{sereq3}
\sum_{n=0}^{\infty} \frac{(q;q^2)_{2n}q^{2n^2+4n}(-1)^n}
{(q^8;q^8)_{n}(-q^2;q^4)_{n+1}} = (-q^9,-q^7,q^8;q^8)_{\infty}\frac{
(q^2;q^4)_{\infty}}{ (q^4;q^4)_{\infty}}.
\end{equation}
\begin{equation}\label{sereq4}
\sum_{n=0}^{\infty} \frac{(q;q^2)_{2n}q^{2n^2}(-1)^n}
{(q^8;q^8)_{n}(-q^2;q^4)_{n}} = (-q^3,-q^5,q^8;q^8)_{\infty}\frac{
(q^2;q^4)_{\infty}}{ (q^4;q^4)_{\infty}}.
\end{equation}
\end{theorem}

\begin{proof} For \eqref{sereq1},
 replace $q$ by $q^2$ in \eqref{eq2}, and set $a=-q^2$ and $c=q^3$.
 Then
 {\allowdisplaybreaks
\begin{align*}
&\sum_{n=0}^{\infty} \frac{(-q^2;q^2)_{n}q^{n^2}}{(q;q)_{2n+1}}=
\frac{1}{1-q}
\sum_{n=0}^{\infty}\frac{(-q^2;q^2)_{n}q^{n^2}}{(q^3;q^2)_{n}(q^2;q^2)_{n}}\\
&=\frac{1}{1-q}\frac{(-q;q^2)_{\infty}} {(q^3;q^2)_{\infty}} =
\frac{(-q;q^2)_{\infty}} {(q;q^2)_{\infty}}.
\end{align*} }

To prove \eqref{sereq2}, we  use the following identity ((25) from
Slater's list, with $q$ replaced by $-q$):
\begin{equation}\label{a25}
\sum_{n=0}^{\infty}\frac{(-1)^n q^{n^2}(q;q^2)_{n}}{(q^4;q^4)_{n}}
=\frac{(-q^3;q^6)_{\infty}^2(q^6;q^6)_{\infty}(q;q^2)_{\infty}}{(q^2;q^2)_{\infty}}.
\end{equation}
In \eqref{eq1}, replace $q$ by $q^2$ and set $a=-1$, $b=q$ and
$\gamma = q^2$. Then
 {\allowdisplaybreaks
\begin{align*}
&\sum_{n=0}^{\infty}
\frac{(-1;q)_{2n}q^{n^2+n}}{(q^2;q^2)_{n}(q^2;q^4)_{n}}\\&=
\sum_{n=0}^{\infty}
\frac{(-1;q^2)_{n}(-q;q^2)_{n}q^{n^2+n}}{(q^2;q^2)_{n}(q;q^2)_{n}(-q;q^2)_{n}}\\&
= \sum_{n=0}^{\infty}
\frac{(-1;q^2)_{n}q^{n^2-n}(q^2)^n}{(q;q^2)_{n}(q^2;q^2)_{n}}\\
&=\frac{(-q^2;q^2)_{\infty}}{(q;q^2)_{\infty}}\sum_{n=0}^{\infty}\frac{(-q)^n
q^{n^2-n}(q;q^2)_{n}}{(-q^2;q^2)_{n}(q^2;q^2)_{n}}
\\
&=\frac{(-q^2;q^2)_{\infty}}{(q;q^2)_{\infty}}\sum_{n=0}^{\infty}\frac{(-1)^n
q^{n^2}(q;q^2)_{n}}{(-q^2;q^2)_{n}(q^2;q^2)_{n}}
\\&
=\frac{(-q^2;q^2)_{\infty}}{(q;q^2)_{\infty}}\sum_{n=0}^{\infty}\frac{(-1)^n
q^{n^2}(q;q^2)_{n}}{(q^4;q^4)_{n}}\\
&=\frac{(-q^2;q^2)_{\infty}}
{(q;q^2)_{\infty}}\frac{(-q^3;q^6)_{\infty}^2(q^6;q^6)_{\infty}(q;q^2)_{\infty}}{(q^2;q^2)_{\infty}}
 \,\,\, (\text{by \eqref{a25}}).
\end{align*} }
The result at \eqref{sereq2} now follows.

Before proving \eqref{sereq3} and \eqref{sereq3}, we recall two
identities from \cite{S52} (identities (38) and (39) (the latter was
also stated by Jackson \cite{J28}):
\begin{equation}\label{S38}
\sum_{n=0}^{\infty}\frac{q^{2n(n+1)}}{(q;q)_{2n+1}}=
\frac{(q^3,q^5,q^{8};q^8)_{\infty}(q^2,q^{14};q^{16})_{\infty}}{(q;q)_{\infty}},
\end{equation}
\begin{equation}\label{S39}
\sum_{n=0}^{\infty}\frac{q^{2n^2}}{(q;q)_{2n}}=
\frac{(q,q^7,q^{8};q^8)_{\infty}(q^6,q^{10};q^{16})_{\infty}}{(q;q)_{\infty}}.
\end{equation}
To prove \eqref{sereq3}, we replace $q$ by $q^2$ in \eqref{eq16} and
then let $c=-q^4$ and $d=q$ to get
\begin{align*}
\sum_{n=0}^{\infty} \frac{(q;q^2)_{2n}q^{2n^2+4n}(-1)^n}
{(q^4;q^4)_{n}(-q^4;q^2)_{2n}}&=
\frac{(q^6;q^4)_{\infty}}{(-q^4;q^2)_{\infty}}\sum_{n=0}^{\infty}
\frac{q^{2n^2+2n}}
{(q^2;q^2)_{n}(q^3;q^2)_{n}}\\
\Longrightarrow \sum_{n=0}^{\infty}
\frac{(q;q^2)_{2n}q^{2n^2+4n}(-1)^n}
{(q^8;q^8)_{n}(-q^2;q^4)_{n+1}}&=
\frac{(q^2;q^4)_{\infty}}{(1+q)(-q^2;q^2)_{\infty}}
\sum_{n=0}^{\infty}\frac{q^{2n^2+2n}} {(q;q)_{2n+1}}\\
&= \frac{(q^2;q^4)_{\infty}}{(1+q)(-q^2;q^2)_{\infty}}
\frac{(q^3,q^5,q^{8};q^8)_{\infty}(q^2,q^{14};q^{16})_{\infty}}{(q;q)_{\infty}},
\end{align*}
where the last equality follows from \eqref{S38}. The result now
follows, after some elementary infinite product manipulations. The
proof of \eqref{sereq4} is similar, except we set $c=-q^2$, $d=q$
and use \eqref{S39}.
\end{proof}

Remark: Shortly after proving \eqref{sereq1}, we discovered, while
reading a pre-print version of \cite{AB07}, that it had previously
been stated by Ramanujan (see \cite[\textbf{Entry 1.7.13}]{AB07}).
However, we include since we discovered it independently and
\cite{AB07} has not yet been published.

\subsection{Watson's Transformation}
Before proving the next identity, we introduce some notation.

An $_{r} \phi _{s}$ \emph{basic hypergeometric series} is defined by
{\allowdisplaybreaks
\begin{multline*} _{r} \phi _{s} \left (
\begin{matrix}
a_{1}, a_{2}, \dots, a_{r}\\
b_{1}, \dots, b_{s}
\end{matrix}
; q,x \right ) := \\
\sum_{n=0}^{\infty} \frac{(a_{1};q)_{n}(a_{2};q)_{n}\dots
(a_{r};q)_{n}} {(q;q)_{n}(b_{1};q)_{n}\dots (b_{s};q)_{n}} \left(
(-1)^{n} q^{n(n-1)/2} \right )^{s+1-r}x^{n},
\end{multline*}
}for $|q|<1$.

Watson's identity is the following. {\allowdisplaybreaks
\begin{multline}\label{Wat}
_{8} \phi _{7} \left (
\begin{matrix}
a,\,q\sqrt{a},\,-q\sqrt{a},\,b,\,c,\,d,\,e,\,q^{-n}\\
\sqrt{a},\,-\sqrt{a},\,aq/b,\,aq/c,\,aq/d,\,aq/e,\,aq^{n+1}\,
\end{matrix}
; q,\frac{a^{2}q^{n+2}}{bcde} \right ) = \\
 \frac{(a q)_{n} (a q/de)_{n}  }{(a q/d)_{n} (a q/e)_{n}
}   \,\,   _{4}\phi _{3} \left (
\begin{matrix}
aq/bc,d,e,q^{-n}\\
aq/b,aq/c,deq^{-n}/a
\end{matrix}\,
; q,q \right ),
\end{multline}
} where $n$ is a non-negative integer. Watson \cite{W29} used his
transformation in his proof of the Rogers-Ramanujan identities
\eqref{rr1}.

\begin{theorem}\label{c12}
Let $a$, $b$ and $q \in \mathbb{C}$, with $|q|<1$. Then
\begin{multline}\label{c12eq}
\sum_{r=0}^{\infty}
\frac{(1+aq^{r})(a^2;q)_{r}(b;q)_{r}\left(-a/b\right)^rq^{r(r+1)/2}}
{(a^2q/b;q)_{r}(q;q)_{r}}\\ =
\frac{(-a;q)_{\infty}(a^2q;q^2)_{\infty}(aq/b;q)_{\infty}}
{(a^2q/b;q)_{\infty}}.
\end{multline}
\end{theorem}
\begin{proof}
We will use Bailey's identity \cite{B41}:
\begin{equation}\label{beq}
\sum_{n=0}^{\infty}
\frac{(a;q)_{n}(b;q)_{n}}{(aq/b;q)_{n}(q;q)_{n}}\left(-\frac{q}{b}\right)^{n}
= \frac{(aq;q^2)_{\infty}(-q;q)_{\infty}(aq^2/b^2;q^2)_{\infty}}
{(aq/b;q)_{\infty}(-q/b;q)_{\infty}}.
\end{equation}
First, let $n, b \to \infty$ in \eqref{Wat} to get
\begin{multline}\label{W3}
\sum_{r \geq 0} \frac{(1-aq^{2r})(a)_{r}(c)_{r}(d)_{r}(e)_{r}\left(
a^{2}/cde \right )^{r} q^{r(r-1)+2r}}
{(1-a)(aq/c)_{r}(aq/d)_{r}(aq/e)_{r}(q)_{r}}\\
=
\frac{(aq)_{\infty}(aq/de)_{\infty}}{(aq/d)_{\infty}(aq/e)_{\infty}}
\sum_{r \geq 0}\frac{  (d)_{r}(e)_{r}
(aq/de)^{r}}{(aq/c)_{r}(q)_{r}}.
\end{multline}
Next, replace $a$ by $-a$, set $c=-b$, $d=a$ and $e=b$, so that
 \eqref{W3} becomes
 \begin{multline*}
\sum_{r \geq 0}
\frac{(1+aq^{2r})(-a)_{r}(-b)_{r}(a)_{r}(b)_{r}\left( -a/b^2 \right
)^{r} q^{r(r-1)+2r}}
{(1+a)(aq/b)_{r}(-q)_{r}(-aq/b)_{r}(q)_{r}}\\
=
\frac{(-aq)_{\infty}(-q/b)_{\infty}}{(-q)_{\infty}(-aq/b)_{\infty}}
\sum_{r \geq 0}\frac{  (a)_{r}(b)_{r}
(-q/b)^{r}}{(aq/b)_{r}(q)_{r}},
 \end{multline*}
 and the result now follows from \eqref{beq}, after replacing $b^2$ by $b$
 and $q^2$ by $q$.
 Note that the result holds initially for $|q/b|<1$, and then
 follows for general $b$ by analytic continuation.
\end{proof}

\subsection{Bailey Pairs}\label{BPS}

A pair of sequences $(\alpha_n, \beta_n)$ that satisfy $\alpha_0=1$
and
\begin{equation}\label{bpeq}
\beta_n = \sum_{r=0}^{n} \frac{\alpha_r}{(q;q)_{n-r}(aq;q)_{n+r}}
\end{equation}
is termed a \emph{Bailey pair relative to $a$}. Bailey \cite{B47,
B49} showed that, for such a pair,
\begin{equation}\label{Baileyeq}
\sum_{n=0}^{\infty} (y,z;q)_{n}\left ( \frac{aq}{yz}\right )^{n}
\beta_n = \frac{(aq/y,aq /z;q)_{\infty}}{ (aq, aq/yz;q)_{\infty}}
\sum_{n=0}^{\infty} \frac{(y,z;q)_{n}}{(aq/y,aq/z;q)_n}\left (
\frac{aq}{yz}\right )^{n} \alpha_n.
\end{equation}
We note two special cases which will be needed later. Firstly, upon
letting $y$, $z \to \infty$ we get that
\begin{equation}\label{Baileyeqspeccase}
\sum_{n=0}^{\infty}a^n q^{n^2} \beta_n = \frac{1}{ (aq;q)_{\infty}}
\sum_{n=0}^{\infty} a^n q^{n^2} \alpha_n.
\end{equation}
Secondly, upon setting $y=q^{1/2}$ and letting $z \to \infty$ we get
that
\begin{equation}\label{Baileyeqspeccasen}
\sum_{n=0}^{\infty}(q^{1/2};q)_n (-1)^n a^n q^{n^2/2} \beta_n =
\frac{(a q^{1/2};q)_{\infty}}{ (aq;q)_{\infty}}
\sum_{n=0}^{\infty}\frac{(q^{1/2};q)_n}{(aq^{1/2};q)_n}\, a^n (-1)^n
q^{n^2/2} \alpha_n.
\end{equation}

\begin{lemma}\label{bplem} The pair $(\alpha_n, \beta_n)$ is a Bailey pair
relative to 1, where
\begin{align*}
\alpha_n &=
\begin{cases} 1, &n=0,\\
2(-1)^n q^{n^2/2},&n\geq 1,
\end{cases}\\
\beta_n &= \frac{(\sqrt{q};q)_n}{(-\sqrt{q},-q,q;q)_n}.
\end{align*}
\end{lemma}
\begin{proof}
Set $a=1$, $c=-\sqrt{q}$, $d=-1$ in Slater's equation (4.1) from
\cite[page 468]{S51}:
\begin{equation*}
\sum_{r=0}^{n} \frac{(1-a q^{2r})(a,c,d;q)_r
q^{(r^2+r)/2}}{(a;q)_{n+r+1}(q;q)_{n-r}(aq/c,aq/d,q;q)_r}\left(\frac{-a}{cd}\right)^r
=\frac{(aq/cd;q)_n}{(aq/c,aq/d,q;q)_n}.
\end{equation*}
The result follows from \eqref{bpeq}, after a little simplification.
\end{proof}

\begin{theorem}\label{c13}
\begin{equation}\label{c13eq}
\sum_{n=0}^{\infty}\frac{q^{2n^2}(q;q^2)_{n}}{(-q;q^2)_{n}(q^4;q^4)_{n}}
=\frac{(q^3;q^6)_{\infty}^{2}(q^6;q^6)_{\infty}}
{(q^2;q^2)_{\infty}}.
\end{equation}
\end{theorem}

\begin{proof}
Substitute the Bailey pair from Lemma \ref{bplem} into
\eqref{Baileyeqspeccase}, with $a=1$, and replace $q$ with $q^2$.
The result follows after using using the Jacobi triple product
identity
\begin{equation}\label{JTP}
\sum_{n=-\infty}^{\infty} z^nq^{n^2}=(-q/z,-qz,q^2;q^2)_{\infty}
\end{equation}
to sum the resulting right side.
\end{proof}

Remark: This is a companion identity to number (27) on Slater's
list, with $q$ replaced by $-q$:
\begin{equation*}
\sum_{n=0}^{\infty}\frac{q^{2n^2+2n}(q;q^2)_{n}}{(-q;q^2)_{n}(q^4;q^4)_{n}}
=\frac{(q;q^6)_{\infty}(q^5;q^6)_{\infty}(q^6;q^6)_{\infty}}
{(q^2;q^2)_{\infty}}.
\end{equation*}

\subsection{An identity of Bailey}
Before coming to the next identity we recall a result of
Bailey(\cite{B51}, p. 220):
\begin{equation}\label{baileyeq}
(-z^2q,-z^{-2}q^3,q^4;q^{4})_{\infty}
+z(-z^2q^3,-z^{-2}q,q^4;q^{4})_{\infty}=(-z,-z^{-1}q,q;q)_{\infty}.
\end{equation}

We also recall Slater's Bailey pair \textbf{G3} (relative to 1) from
\cite{S51}.
\begin{align}\label{bpg3}
\alpha_n&=
\begin{cases}
1,&n=0,\\
q^{3r^2}(q^{3r/2}+q^{-3r/2}),&n=2r,\,\, r\geq1,\\
-q^{3r^2}(q^{3r/2}+q^{9r/2+3/2}),&n=2r+1,
\end{cases}\\
\beta_n&=\frac{q^{n}}{(q^2;q^2)_n(-q^{1/2};q)_n}. \notag
\end{align}
We note that Slater used \textbf{G3} to derive two other
series-product identities, (16) and (32) in \cite{S52}, so we may
regard the identity in Theorem \ref{tqd} as one she missed.
\footnote{ In an earlier version of this paper we proved Theorem
\ref{tqd} by the method of $q$-difference equations. However, that
proof was much longer and less transparent than the present proof.}

\begin{theorem}\label{tqd}
Let $|q|<1$. Then
\begin{equation}\label{qdeq}
\sum_{n=0}^{\infty}\frac{(-1)^{n}q^{n^2+2n}(q;q^2)_{n}}{(-q;q^2)_{n}(q^4;q^4)_{n}}
=\frac{(-q;q^5)_{\infty}(-q^4;q^5)_{\infty}(q^5;q^5)_{\infty}(q;q^2)_{\infty}}
{(q^2;q^2)_{\infty}}.
\end{equation}
\end{theorem}
Remark: This identity is clearly a companion to Identity (21) on
Slater's list:
\begin{equation}\label{id21eq}
\sum_{n=0}^{\infty}\frac{(-1)^{n}q^{n^2}(q;q^2)_{n}}{(-q;q^2)_{n}(q^4;q^4)_{n}}
=\frac{(-q^3;q^5)_{\infty}(-q^2;q^5)_{\infty}(q^5;q^5)_{\infty}(q;q^2)_{\infty}}
{(q^2;q^2)_{\infty}}.
\end{equation}

\begin{proof}[Proof of Theorem \ref{tqd}]
We insert the Bailey pair \eqref{bpg3} into
\eqref{Baileyeqspeccasen}, set $a=1$ and replace $q$ by $q^2$ to get
\begin{align*}
\sum_{n=0}^{\infty}&\frac{(-1)^{n}q^{n^2+2n}(q;q^2)_{n}}{(-q;q^2)_{n}(q^4;q^4)_{n}}\\
&=\frac{(q;q^2)_{\infty}}{(q^2;q^2)_{\infty}} \left (1+
\sum_{r=1}^{\infty}q^{10r^2}(q^{3r}+q^{-3r})+\sum_{r=0}^{\infty}q^{10r^2+4r+1}(q^{3r}+q^{9r+3})\right)\\
&=\frac{(q;q^2)_{\infty}}{(q^2;q^2)_{\infty}}
\left(\sum_{r=-\infty}^{\infty}q^{10r^2+3r}+q^4\sum_{r=-\infty}^{\infty}q^{10r^2+13r}\right)\\
&=\frac{(q;q^2)_{\infty}}{(q^2;q^2)_{\infty}}
\left((-q^{7},-q^{13},q^{20};q^{20})_{\infty}+q^{4}(-q^{-3},-q^{23},q^{20};q^{20})_{\infty}\right)\\
&=\frac{(-q^3;q^5)_{\infty}(-q^2;q^5)_{\infty}(q^5;q^5)_{\infty}(q;q^2)_{\infty}}
{(q^2;q^2)_{\infty}}.
\end{align*}
The next-to-last equation follows from \eqref{JTP}, and the last
equation  follows from \eqref{baileyeq}, with $q$ replaced by $q^5$
and $z=q^4$.
\end{proof}

\subsection{Miscellaneous Methods}
Before coming to the next identity, we recall two other necessary
results. The first of these is an identity of Blecksmith, Brillhart
and Gerst \cite{BBG88} (a proof is also given in \cite{CH01}):
\begin{equation}\label{bbgid}
\sum_{n=-\infty}^{\infty} q^{n^2}-\sum_{n=-\infty}^{\infty} q^{5n^2}
=2 q \frac{(q^4,q^6,q^{10},q^{14},q^{16},q^{20};q^{20})_{\infty}}
{(q^3,q^7,q^{8},q^{12},q^{13},q^{17};q^{20})_{\infty}}.
\end{equation}
If we replace $q$ by $-q$ and apply the Jacobi triple product
identity to the left side, \eqref{bbgid} may be re-written as
{\allowdisplaybreaks
\begin{equation}\label{bbgid2}
(q^5, q^5,q^{10};q^{10})_{\infty}-(q,q,q^{2};q^{2})_{\infty} =2 q
\frac{(q^4,q^6,q^{10};q^{10})_{\infty}} {(-q^3,-q^7;q^{10})_{\infty}
(q^{8},q^{12};q^{20})_{\infty}}.
\end{equation}
}

The second is the following identity, due to Rogers~\cite[p. 330
(4), line 3, corrected]{R17} recently generalized by the third
author~\cite[p. 404, Eq. (3)]{S06}:
\begin{equation}\label{tseq}
\sum_{n=0}^{\infty}\frac{q^{n(n+1)/2}(-1;q)_{n}}{(q;q)_{n}(q;q^2)_{n}}
=\frac{(q^5,q^5,q^{10};q^{10})_{\infty}}
{(q;q)_{\infty}(q;q^2)_{\infty}}.
\end{equation}

We are now able to prove another identity discovered during the
present investigations.
\begin{theorem}\label{tq}Let $|q|<1$. Then
{\allowdisplaybreaks
\begin{equation}\label{c15eq}
\sum_{n=0}^{\infty}\frac{q^{(n^2+3n)/2}(-q;q)_{n}}
{(q;q^2)_{n+1}(q;q)_{n+1}} =\frac{(q^{10};q^{10})_{\infty}}
{(q;q)_{\infty}(q;q^2)_{\infty}
(-q^3,-q^4,-q^6,-q^7;q^{10})_{\infty}}.
\end{equation}
}
\end{theorem}

\begin{proof}
{\allowdisplaybreaks
\begin{align*}
\sum_{n=0}^{\infty}&\frac{q^{(n^2+3n)/2}(-q;q)_{n}}
{(q;q^2)_{n+1}(q;q)_{n+1}} = \frac{1}{2q}
\sum_{n=0}^{\infty}\frac{q^{(n+1)(n+2)/2}(-1;q)_{n+1}}
{(q;q^2)_{n+1}(q;q)_{n+1}}\\
&= \frac{1}{2q}\left (
\sum_{n=0}^{\infty}\frac{q^{n(n+1)/2}(-1;q)_{n}}
{(q;q^2)_{n}(q;q)_{n}}-1 \right )\\
&= \frac{1}{2q}\left (\frac{(q^5,q^5,q^{10};q^{10})_{\infty}}
{(q;q)_{\infty}(q;q^2)_{\infty}} -1 \right ) \, (\text{by
\eqref{tseq}} )\\
&= \frac{1}{2q(q;q)_{\infty}(q;q^2)_{\infty}}\left
((q^5,q^5,q^{10};q^{10})_{\infty} -(q;q^2)_{\infty}
(q;q^2)_{\infty}(q^2;q^2)_{\infty}\right ) \\
&=\frac{(q^4,q^6,q^{10};q^{10})_{\infty}}
{(q;q)_{\infty}(q;q^2)_{\infty}(-q^3,-q^7;q^{10})_{\infty}
(q^{8},q^{12};q^{20})_{\infty}}\, (\text{by \eqref{bbgid2}})\\
&=\frac{(q^{10};q^{10})_{\infty}}
{(q;q)_{\infty}(q;q^2)_{\infty}(-q^3,-q^4,-q^6,-q^7;q^{10})_{\infty}
}.
\end{align*}
}
\end{proof}

\section{An Application of the Two-Variable Generalization of
a Rogers-Ramanujan Type Series}
  In~\cite{A86}, Andrews showed that a certain two-variable generalization
$f(t,q)$ of a Rogers-Ramanujan type series $\Sigma (q)$ served as a generating
function in $t$ of a sequence of polynomials $P_n(q)$ for which
$\lim_{n\to\infty} P_n(q) = \Sigma (q)$.

   For example, if
   \begin{equation} \label{RR1ser}
  \Sigma:=\Sigma(q):= \sum_{n=0}^\infty \frac{q^{n^2}}{(q;q)_n},
 \end{equation}
 the series associated with the first Rogers-Ramanujan identity,
and the two variable generalization $f(t,q)$ of $\Sigma$ is given by
\begin{equation} \label{tRR1ser}
  f(t,q):= \sum_{n=0}^\infty \frac{ t^{2n} q^{n^2} }{ (t;q)_{n+1}},
\end{equation}
then it is also the case that
\begin{equation}
  f(t,q) = \sum_{n=0}^\infty P_n (q) t^n
\end{equation}
where
\begin{equation} \label{SchurPoly} P_0 (q) = P_1 (q) = 1; \qquad
P_n(q) = P_{n-1} (q) + q^{n-1} P_{n-2} (q) \quad \mbox{if $n\geq 2$}.
\end{equation}
Note that the polynomials in~\eqref{SchurPoly}, which are
$q$-analogs of the Fibonacci numbers, are sometimes
called the \emph{Schur polynomials} because they were employed
by Schur in his proof of the Rogers-Ramanujan identities.
Indeed, Schur~\cite{S17} showed that
 \begin{equation} \label{SchurRep}
 P_n(q) = \sum_{j=-\infty}^\infty (-1)^j q^{j(5j+1)/2}
 \qbin{n}{\lfloor  \frac{n+5j+1}{2}  \rfloor}{q},
 \end{equation}
 while elsewhere MacMahon~\cite{M18} showed that
 \begin{equation}\label{PAMRep}
   P_n(q) = \sum_{j\geq 0} q^{j^2} \qbin{n-j}{j}{q},
 \end{equation}
 where
 \[ \qbin{A}{B}{q} := \left\{
   \begin{array}{ll}
      (q;q)_A (q;q)_{B}^{-1} (q;q)_{A-B}^{-1} &\mbox{if $0\leq B\leq A$} \\
      0 & \mbox{otherwise}. \end{array} \right.
       \]
 Thus, by combining~\eqref{PAMRep} and~\eqref{SchurRep}, we
 may observe, as Andrews did in~\cite{A70}, that we have
 a polynomial identity which generalizes the first Rogers-Ramanujan
 identity and that we may recover~\eqref{rr1} by letting
 $n\to\infty$.

    In~\cite{S03}, the third author used these two-variable generalizations
  $f(t,q)$
 of Rogers-Ramanujan type series to find polynomial generalizations
 of all 130 identities in Slater's list~\cite{S52}.  Previously, Santos~\cite{S91}
 had studied a large number of polynomial sequences associated with
 two-variable generalizations of series in Slater's list.  Indeed, the primary
 use of the $f(t,q)$ has been as a generating function in $t$ for
 sequences of polynomials.

 However, here we wish to turn our attention to a different use of
the $f(t,q)$ by following up on an observation made by
Andrews~\cite[p. 89]{A86}.  Letting $f(q,t)$ and $\Sigma(q)$ be
as above, Andrews noted the following:
not only do we have, as required,
  \begin{equation*}
   \lim_{t\to 1^{-}} (1-t) f(t,q) = \Sigma(q),
  \end{equation*}
but also
   \begin{equation*}
     \lim_{t\to -1^+} (1-t) f(t,q) = f_0 (q),
   \end{equation*}
where
  \begin{equation*}
  f_0(q) := \sum_{n=0}^\infty \frac{q^{n^2}}{(-q;q)_n},
   \end{equation*}
one of Ramanujan's fifth order mock theta functions (cf.~\cite{W36},
~\cite{W37}).

  Elsewhere~\cite[p. 90--91]{A86}, Andrews notes that if we take
$\Sigma$ to the be the Rogers-Ramanujan type series associated with
Eq. (46) on Slater's list~\cite{S52}, and define its two variable generalization
as
\begin{equation*}
  f(t,q):= \sum_{n=0}^\infty \frac{ t^{3n} q^{n(3n-1)/2}}{ (t;q)_{n+1}
(t^2q;q^2)_n},
\end{equation*}
we find that $\lim_{t\to -1^+} f(t,q)$ is not a mock theta function,
but rather a~\emph{false} theta function
studied by Rogers~\cite[p. 333(2)]{R17}.

  Andrews later comments~\cite[p. 93]{A86}: ``Now if we view this
as a curve $y=f(t)$ the points of which are functions of $q$, we
find that frequently if $f(1)$ is a modular form, $f(-1)$ is a mock
or false theta function.  Is there some general structure possible
in which this seemingly amazing occurrence becomes more explicable?"
While we do not have an answer to Andrews' question, we have
observed that there is a third possibility.  Namely, that $f(1)$ is
a modular form and $f(-1)$ neither a mock nor false theta function,
but rather a sum of modular forms.

  Consider the following family of Rogers-Ramanujan type identities
related to the modulus 27 which are due to
Dyson~\cite[p. 433, Eqs. (B1)--(B4)]{B47}
and reproved by Slater~\cite[p. 161--2, Eqs. (90)--(93)]{S52}.
\begin{align}
1+ \sum_{n=1}^\infty \frac{ q^{n^2} (q^3;q^3)_{n-1} }{(q;q)_n (q;q)_{2n-1}}
 &= \frac{(q^{12}, q^{15}, q^{27}; q^{27})_\infty}{(q;q)_\infty} \label{Dyson4}\\
   \sum_{n=0}^\infty \frac{ q^{n(n+1)} (q^3;q^3)_n }{(q;q)_n (q;q)_{2n+1}}
 &= \frac{(q^9 ; q^{9})_\infty}{(q;q)_\infty} \label{Dyson3} \\
  \sum_{n=0}^\infty \frac{ q^{n(n+2)} (q^3;q^3)_n }{(q;q)_n (q;q)_{2n+2}}
 &= \frac{(q^6, q^{21}, q^{27}; q^{27})_\infty}{(q;q)_\infty} \label{Dyson2}\\
  \sum_{n=0}^\infty \frac{  q^{n(n+3)} (q^3;q^3)_n }{(q;q)_n (q;q)_{2n+2}}
 &= \frac{(q^3, q^{24}, q^{27}; q^{27})_\infty}{(q;q)_\infty} \label{Dyson1}
\end{align}

The relevant two-variable generalizations
(see~\cite[p. 15, Thm. 2.2]{S03}) of~\eqref{Dyson4}--\eqref{Dyson1}
are
\begin{align}
f_{\ref{Dyson4}}(t,q) &:=
\frac{1}{1-t} +\sum_{n=1}^\infty \frac{ t^{2n} q^{n^2}
(t^3 q^3;q^3)_{n-1} }{(t;q)_{n+1} (t^2 q;q)_{2n-1}}\label{fDyson4}\\
f_{\ref{Dyson3}}(t,q) &:=  \sum_{n=0}^\infty \frac{ t^{2n} q^{n(n+1)}
(t^3 q^3;q^3)_n }{(t;q)_{n+1} (t^2 q;q)_{2n+1}}
  \label{fDyson3}\\
   f_{\ref{Dyson2}}(t,q) &:=  \sum_{n=0}^\infty \frac{ t^{2n} q^{n(n+2)}
(t^3 q^3;q^3)_n }{(t;q)_{n+1} (t^2 q;q)_{2n+2}}
  \label{fDyson2}\\
 f_{\ref{Dyson1}}(t,q) &:=  \sum_{n=0}^\infty \frac{ t^{2n} q^{n(n+3)}
(t^3 q^3;q^3)_n }{(t;q)_{n+1} (t^2 q;q)_{2n+2}}
  \label{fDyson1}.
\end{align}

We believe that the four identities related to the modulus 108 recorded
below, which arise from the $t=-1$ cases
of~\eqref{fDyson4}--\eqref{fDyson1}, are new.

\begin{multline}
 1+\sum_{n=1}^\infty \frac{ q^{n^2} (-q^3;q^3)_{n-1} }{(-q;q)_n (q;q)_{2n-1}}
 =\\
\frac{(q^{12},q^{15},q^{27};q^{27})_\infty-2q^2
(-q^{33},-q^{75},q^{108};q^{108})_\infty+
2q^{7}(-q^{15},-q^{93},q^{108};q^{108})_\infty}{(q;q)_\infty}
\label{CDyson4}
\end{multline}
\begin{multline}
 \sum_{n=0}^\infty \frac{ q^{n(n+1)} (-q^3;q^3)_{n} }{(-q;q)_n (q;q)_{2n+1}}
  =\\\frac{(q^9,q^{18},q^{27};q^{27})_\infty -2q^3
(-q^{27},-q^{81},q^{108};q^{108})_\infty + 2q^{9}
 (-q^{9},-q^{99},q^{108};q^{108})_\infty}{(q;q)_\infty}
\label{CDyson3}
\end{multline}
\begin{multline}
 \sum_{n=0}^\infty \frac{ q^{n(n+2)} (-q^3;q^3)_n }{(-q;q)_n (q;q)_{2n+2}}
=\\
\frac{(q^6,q^{21},q^{27};q^{27})_\infty -2q^4
(-q^{21},-q^{87},q^{108};q^{108})_\infty + 2q^{11}
 (-q^{3},-q^{105},q^{108};q^{108})_\infty}{(q;q)_\infty}
\label{CDyson2}
\end{multline}
\begin{multline}
\sum_{n=0}^\infty \frac{  q^{n(n+3)} (-q^3;q^3)_n }{(-q;q)_n
(q;q)_{2n+2}}
\\
=\text{\small{$\frac{(q^3,q^{24},q^{27};q^{27})_\infty-2q^5(-q^{15},-q^{93},q^{108};q^{108})_\infty+2q^{13}
(-q^{-3},-q^{111}, q^{108};q^{108})_\infty}{(q;q)_\infty}
\label{CDyson1}$}}
\end{multline}

The identities~\eqref{CDyson4}--\eqref{CDyson1} may be proved using
Bailey pairs (see Sec.~\ref{BPS}).
Although less well known than~\eqref{bpeq}, the following characterization of
Bailey pairs~\cite[p. 29, Eq. (3.40) with $a=1$]{A86}
is equivalent to~\eqref{bpeq} with $a=1$:
\begin{equation} \label{BPequiv}
\alpha_n = (1-q^{2n}) \sum_{k=0}^n
  \frac{(-1)^{n-k} q^{\binom{n-k}{2}} (q;q)_{n+k-1}} {(q;q)_{n-k} } \beta_k
\end{equation}
Furthermore, it is straightforward to show that~\eqref{BPequiv} may be
rewritten as
\begin{equation} \label{BPdef2}
\alpha_n = (-1)^n q^{\binom{n}{2}} (1+q^n) \sum_{k=0}^n
(q^n;q)_k (q^{-n};q)_k q^k \beta_k,
\end{equation}
which is the form we shall employ.

\begin{lemma} \label{BaileyPair}
If, for $n$ a nonnegative integer,
\begin{equation} \label{alpha}
  \alpha_n = \left\{
    \begin{array}{ll}
       1 &\mbox{if $n=0$}\\
       (-1)^r q^{\frac 92 r^2 - \frac 32 r} (1+q^{3r}) &\mbox{if $n=3r>0$}\\
      -2q^{18r^2+9r+1} &\mbox{if $n=6r+1$}\\
      2 q^{18r^2+15r+3} &\mbox{if $n=6r+2$}\\
      2 q^{18r^2+ 21r+6} &\mbox{if $n=6r+4$}\\
      -2 q^{18r^2+27r+10} &\mbox{if $n=6r+5$}
    \end{array} \right.
\end{equation} and
\begin{equation} \label{beta}
  \beta_n = \left\{
    \begin{array}{ll}
       \frac{(-q^3;q^3)_{n-1}}{(-q;q)_n (q;q)_{2n-1} } &\mbox{if $n>0$}\\
       1 &\mbox{if $n=0$,}
    \end{array} \right.
\end{equation} then $( \alpha_n, \beta_n )$ forms a Bailey pair.
\end{lemma}
\begin{proof}
   By inserting~\eqref{alpha} and~\eqref{beta} into~\eqref{BPdef2},
 it is clear that we will be done once we show
    \begin{multline} \label{first}
    (-1)^r q^{\frac 92 r^2 - \frac 32r} (1+q^{3r})
      \\= (-1)^r q^{\binom{3r}{2}} (1+q^{3r} )\left\{ 1+ \sum_{k=1}^{3r}
\frac{(q^{3r};q)_k (q^{-3r};q)_k  (-q^3;q^3)_k q^k}{(-q;q)_k (q;q)_{2k-1} } \right\},
\end{multline}
\begin{multline}
  -2q^{18r^2+9r+1} \\
  =  -q^{\binom{6r+1}{2}} (1+q^{6r+1} )\left\{ 1+ \sum_{k=1}^{6r+1},
\frac{(q^{6r+1};q)_k (q^{-6r-1};q)_k  (-q^3;q^3)_k q^k}{(-q;q)_k (q;q)_{2k-1} } \right\},
\end{multline}
\begin{multline}
    2 q^{18r^2+15r+3} \\=  q^{\binom{6r+2}{2}} (1+q^{6r+2} )\left\{ 1+ \sum_{k=1}^{6r+2}
\frac{(q^{6r+2};q)_k (q^{-6r-2};q)_k  (-q^3;q^3)_k q^k}{(-q;q)_k (q;q)_{2k-1} } \right\} ,
\end{multline}
\begin{multline}
    2 q^{18r^2+ 21r+6}\\=  q^{\binom{6r+4}{2}} (1+q^{6r+4} )\left\{ 1+ \sum_{k=1}^{6r+4}
\frac{(q^{6r+4};q)_k (q^{-6r-4};q)_k  (-q^3;q^3)_k q^k}{(-q;q)_k (q;q)_{2k-1} } \right\},
\end{multline}
and
\begin{multline}
    -2 q^{18r^2+27r+10}\\ =  -q^{\binom{6r+5}{2}} (1+q^{6r+5} )\left\{ 1+ \sum_{k=1}^{6r+5}
\frac{(q^{6r+5};q)_k (q^{-6r-5};q)_k  (-q^3;q^3)_k q^k}{(-q;q)_k (q;q)_{2k-1} } \right\},
\label{last}
\end{multline}
for nonnegative integers $r$.
Using elementary algebra, the equations~\eqref{first}--\eqref{last} are easily shown
to be equivalent to
\begin{align}
\sum_{k=1}^{3r}
\frac{(q^{3r};q)_k (q^{-3r};q)_k  (-q^3;q^3)_k q^k}{(-q;q)_k (q;q)_{2k-1} }  &=0, \label{eqfirst}\\
\frac{q^{6r+1}+1}{q^{6r+1}-1} \sum_{k=1}^{6r+1}
\frac{(q^{6r+1};q)_k (q^{-6r-1};q)_k  (-q^3;q^3)_k q^k}{(-q;q)_k (q;q)_{2k-1} } &=1,\\
\frac{q^{6r+2}+1}{q^{6r+2}-1} \sum_{k=1}^{6r+2}
\frac{(q^{6r+2};q)_k (q^{-6r-2};q)_k  (-q^3;q^3)_k q^k}{(-q;q)_k (q;q)_{2k-1} } &=1,\\
\frac{1+q^{6r+4}}{1-q^{6r+4}} \sum_{k=1}^{6r+4}
\frac{(q^{6r+4};q)_k (q^{-6r-4};q)_k  (-q^3;q^3)_k q^k}{(-q;q)_k (q;q)_{2k-1} } &=1,\\
\frac{1+q^{6r+5}}{1-q^{6r+5}} \sum_{k=1}^{6r+5}
\frac{(q^{6r+5};q)_k (q^{-6r-5};q)_k  (-q^3;q^3)_k q^k}{(-q;q)_k (q;q)_{2k-1} } &=1.
\label{eqlast}
\end{align}

Each of equations~\eqref{eqfirst}--\eqref{eqlast} may be verified using the
WZ method~\cite[Chapter 7]{PWZ96}
or induction on $r$.
\end{proof}

\begin{theorem}
Identity~\eqref{CDyson4} is valid.
\end{theorem}
\begin{proof}
Recall that the weak form of Bailey's lemma~\cite[p. 27, Eq. (3.33) with $a=1$]{A86}
states that
\begin{equation}\label{WBL}
  \sum_{n=0}^\infty q^{n^2} \beta_n = \frac{1}{(q;q)_\infty} \sum_{n=0}^\infty q^{n^2} \alpha_n
\end{equation} for any Bailey pair $( \alpha_n, \beta_n )$.
Inserting the Bailey pair established in Lemma~\ref{BaileyPair} into~\eqref{WBL} yields
\begin{multline}
  1+\sum_{n=1}^\infty \frac{ q^{n^2} (-q^3;q^3)_{n-1}}{(-q;q)_n (q;q)_{2n-1}} \\=
  \frac{1}{(q;q)_\infty} \left(  1 + \sum_{r=1}^\infty (-1)^r q^{\frac{27}{2}r^2- \frac 32 r}(1+q^{3r})
  -2 \sum_{r=0}^\infty  q^{54r^2+21r+2} \right. \\ \left. + 2\sum_{r=0}^\infty q^{54r^2+39r+7}
  +2\sum_{r=0}^\infty q^{54r^2+69r+22} -2\sum_{r=0}^\infty q^{54r^2+87r+35} \right)\\
  =  \frac{1}{(q;q)_\infty} \left( \sum_{r=-\infty}^\infty (-1)^r q^{\frac{27}{2}r^2- \frac 32 r}
  -2 \sum_{r=0}^\infty  q^{54r^2+21r+2}+ 2\sum_{r=0}^\infty q^{54r^2+39r+7}  \right. \\ \left.
  +2\sum_{r=1}^\infty q^{18r^2-39r+7} -2\sum_{r=1}^\infty q^{54r^2-21r+2} \right)\\
  =  \frac{1}{(q;q)_\infty} \left( \sum_{r=-\infty}^\infty (-1)^r q^{\frac{27}{2}r^2- \frac 32 r}
  -2 \sum_{r=-\infty}^\infty  q^{54r^2-21r+2} + 2\sum_{r=-\infty}^\infty q^{54r^2-39r+7}
  \right)\\
   \\\text{\small{$=\frac{(q^{12},q^{15},q^{27};q^{27})_\infty-2q^2(-q^{33},-q^{75},q^{108};q^{108})_\infty
+2q^7(-q^{15},-q^{93};q^{108};q^{108})_\infty }{(q;q)_\infty}.$}}
\end{multline}
\end{proof}
The other identities~\eqref{CDyson3}--\eqref{CDyson1} may be proved similarly.

 The identity~\eqref{CDyson3} deserves special attention because
its right hand side may be expressed as a single infinite product whereas it appears that
of the other three can not be simplified beyond a sum of three infinite products.
\begin{theorem}
\begin{equation}
 \sum_{n=0}^\infty \frac{ q^{n(n+1)} (-q^3;q^3)_{n} }{(-q;q)_n (q;q)_{2n+1}}
 = \frac{(q^3;q^3)_\infty (q^3;q^{18})_\infty (q^{15};q^{18})_\infty}{(q;q)_\infty}
\label{CDyson3-2}
\end{equation}
\end{theorem}
\begin{proof}
We shall require Fricke's quintuple product identity~\cite{C06}
\begin{equation} \label{QPI}
(z^3 q, z^{-3} q^2, q^3;q^3)_\infty + z (z^{-3} q, z^3 q^2,q^3;q^3)_\infty =
  (-z^{-1} q, -z, q; q)_\infty (z^{-2}q, z^2 q; q^2)_\infty
\end{equation} and an identity due to Bailey~\cite[p. 220, Eq. (4.1)]{B51}
\begin{equation}\label{BaileySimp}
  (-z^2 q, -z^{-2} q^3, q^4; q^4)_\infty + z (-z^2 q^3, -z^{-2} q, q^4; q^4)_\infty = (-z, -z^{-1}q,q;q)_\infty.
\end{equation}
By~\eqref{CDyson3}, we have
\begin{multline*}
 \sum_{n=0}^\infty \frac{ q^{n(n+1)} (-q^3;q^3)_{n} }{(-q;q)_n (q;q)_{2n+1}}
 \\ =
 \text{\small{$\frac{(q^9, q^{18}, q^{27}; q^{27})_\infty
          -2q^3 (-q^{27},-q^{81},q^{108};q^{108})_\infty + 2q^{9}
 (-q^{9},-q^{99}, q^{108};q^{108})_\infty}{(q;q)_\infty}.$}}
\label{CDyson3}
\end{multline*}
Expanding the first triple product in the numerator by~\eqref{BaileySimp}
with $q$ replaced by $q^{27}$ and $z=-q^9$ yields
\[\text{\small{$\frac{(-q^{45},-q^{63},q^{108};q^{108})_\infty
-2q^3 (-q^{27},-q^{81},q^{108};q^{108})_\infty
+q^{9}
(-q^{9},-q^{99},q^{108};q^{108})_\infty}{(q;q)_\infty}.$}}\] Two
further applications of~\eqref{BaileySimp} shows that the
preceding expression is equal to
\begin{align*} & \quad \frac{(-q^{9}, -q^{18}, q^{27}; q^{27})_\infty
          -q^3 (-1,-q^{27},q^{27};q^{27})_\infty}{(q;q)_\infty}\\
    &= \frac{(q^3, q^6, q^9;q^9)_\infty (q^3, q^{15}; q^{18})_\infty}{(q;q)_\infty}\\
    & =  \frac{(q^3;q^3)_\infty (q^3, q^{15}; q^{18})_\infty}{(q;q)_\infty} ,
 \end{align*}
 where the penultimate equality follows from~\eqref{QPI}.
\end{proof}

   We close this section by recalling that once a given Bailey pair is established,
it may be utilized in connection with limiting cases of Bailey's
lemma other than \eqref{WBL}, thus yielding additional
Rogers-Ramanujan type identities. For instance, if we insert the
Bailey pair established in Lemma~\ref{BaileyPair} into \cite[p.
26, Eq. (3.28) with $n,\rho_1\to\infty$ and $\rho_2 = -\sqrt{q}$
]{A86}, we obtain the identity related to the modulus 144:
\begin{multline} \label{TNSId}
1+\sum_{n=1}^\infty \frac{ q^{n^2} (-q;q^2)_n (-q^6;q^6)_{n-1} }{ (-q^2;q^2)_n
(q^2;q^2)_{2n-1}} \\
\text{\small{$=\frac{ (q^{15},q^{21},q^{36};q^{36})_\infty - 2q^3
(-q^{42},-q^{102},q^{144};q^{144})_\infty
  +2q^{10} (-q^{18}, -q^{126}, q^{144}; q^{144})_\infty }
{ (q;q^2)_\infty (q^4;q^4)_\infty}.$}}
\end{multline}
A partner of~\eqref{TNSId} is
\begin{multline}
\sum_{n=0}^\infty \frac{ q^{n^2+4n} (-q;q^2)_{n+1} (-q^6;q^6)_{n} }{ (-q^2;q^2)_n
(q^2;q^2)_{2n+2}} \\
\text{\small{$=\frac{ (q^{3},q^{33},q^{36};q^{36})_\infty - 2q^7
(-q^{18},-q^{126},q^{144};q^{144})_\infty
 +2q^{12} (-q^{6}, -q^{138}, q^{144}; q^{144})_\infty }
{ (q;q^2)_\infty (q^4;q^4)_\infty}.$}}
\end{multline}

\section{Concluding Remarks}
For quite a long time we were convinced that there must exist a
general transformation of the type found  in Watson's theorem (see
\eqref{Wat}), a transformation which would give the result in
Theorem \ref{tqd} as a special case for particular values of its
parameters.

One reason we thought this transformation had to exist was the
appearance of the
$(-q;q^5)_{\infty}(-q^4;q^5)_{\infty}(q^5;q^5)_{\infty}$ term on the
product side, which can be represented as an infinite series via the
Jacobi triple product. This in turn brought to mind Watson's proof
of the Rogers-Ramanujan identities, where he showed that these
followed as special cases of \eqref{Wat}.

However, we could not find such a transformation, but possibly our
search was incomplete. Does the identity in Theorem \ref{tqd} follow
as a special case of some known transformation, perhaps some known
transformation between basic hypergeometric series? Is this identity
a special case of some as yet undiscovered general transformation?

As remarked at the end of  the introduction, varying the form of the
series $S$ in \eqref{Seq} may lead to other new identities of the
Rogers-Ramanujan-Slater type. In particular, one might hope for the
discovery of new identities which are not readily proved within the
framework of our present understanding of identities of the
Rogers-Ramanujan-Slater type. We hope to continue these
investigations in a subsequent paper.

 \allowdisplaybreaks{

}
\end{document}